\documentclass[12pt,reqno]{article}

\usepackage[usenames]{color}
\usepackage[colorlinks=true,
linkcolor=webgreen, filecolor=webbrown,
citecolor=webgreen]{hyperref}

\definecolor{webgreen}{rgb}{0,.5,0}
\definecolor{webbrown}{rgb}{.6,0,0}

\usepackage{amssymb}
\usepackage{graphicx}
\usepackage{amscd}
\usepackage{lscape}
\usepackage{tikz}
\usepackage{pgfplots}

\usetikzlibrary{matrix}
\usetikzlibrary{fit,shapes}
\usetikzlibrary{positioning, calc}
\tikzset{circle node/.style = {circle,inner sep=1pt,draw, fill=white},
        X node/.style = {fill=white, inner sep=1pt},
        dot node/.style = {circle, draw, inner sep=5pt}
        }
\usepackage{tkz-fct}

\usepackage{amsthm}
\newtheorem{theorem}{Theorem}

\newtheorem{proposition}[theorem]{Proposition}

\theoremstyle{definition}

\newtheorem{example}[theorem]{Example}

\usepackage{float}

\usepackage{graphics,amsmath}
\usepackage{amsfonts}
\usepackage{latexsym}
\usepackage{epsf}

\DeclareMathOperator{\sech}{sech}
\DeclareMathOperator{\gd}{gd}
\DeclareMathOperator{\erf}{erf}

\begin{document}

\begin{center}
\vskip 1cm{\LARGE\bf Sigmoid functions and exponential Riordan arrays} \vskip 1cm \large
Paul Barry\\
School of Science\\
Waterford Institute of Technology\\
Ireland\\
\href{mailto:pbarry@wit.ie}{\tt pbarry@wit.ie}
\end{center}
\vskip .2 in

\begin{abstract} Sigmoid functions play an important role in many areas of applied mathematics, including machine learning, population dynamics and probability. We place the study of sigmoid functions in the context of the derivative sub-group of the group of exponential Riordan arrays. Links to families of polynomials are drawn, and it is shown that in some cases these polynomials are orthogonal. In the non-orthogonal case, transformations are given that produce orthgonal systems. Alternative means of characterisation are given, based on the production (Stieltjes) matrix associated to the relevant Riordan array. \end{abstract}

\section{Introduction}
Sigmoid functions, familiarly known as ``S'' shaped curves, and  typified by the function
$$f(x)=\frac{1}{1+e^{-x}},$$ occur in many areas of applied mathematics. These areas include machine learning, population dynamics, economics, and probability. Sigmoid functions \cite{Lynch} are assumed to be real-valued and differentiable, and their derivatives exhibit a ``bell-shaped'' curve over the interval of interest. See Figure \ref{sigmoid}.

\begin{center}
\begin{figure}
\begin{center}
\begin{tikzpicture}
\begin{axis} [width=10cm]
\addplot[smooth, very thick, color=black] {1/(1+exp(-x))};
\addplot[smooth, very thick, color=blue] {exp(x)/((1+exp(x))^2};
\end{axis}
\end{tikzpicture}
\end{center}
\caption{The sigmoid pair $(\frac{1}{1+e^{-x}})', \frac{1}{1+e^{-x}}$}\label{sigmoid}
\end{figure}
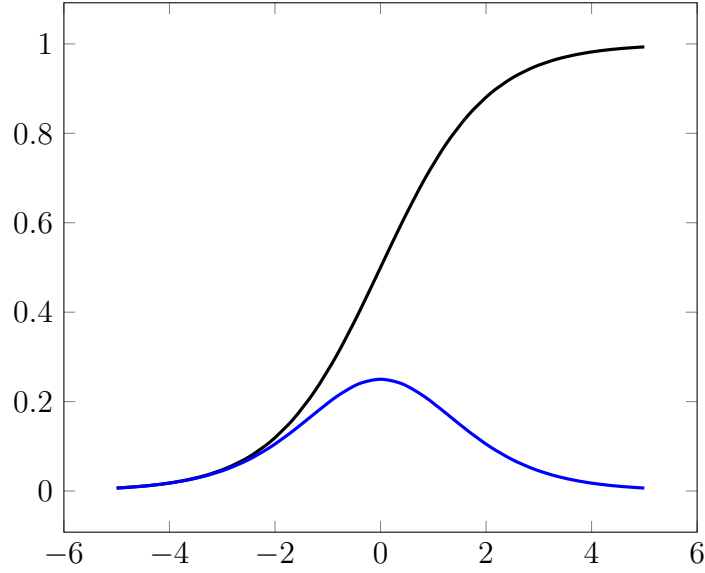
\end{center}

In this note, we shall focus on the pair $(f'(x), f(x))$ where $f(x)$ is a sigmoid function. Examples of choices for $f(x)$ that we shall meet are $\tan^{-1}(x)$, $\tanh(x)$, $\frac{x}{\sqrt{1+x^2}}$, and the Gudermanian function
$\gd(x)=\tan^{-1}(\sinh(x))$.

We shall place the study of such pairs of functions within the context of the derivative sub-group of the group of exponential Riordan arrays \cite{SGWW, Survey}. This will allow us to develop a framework for the investigation of the properties of these pairs of functions, and will help us to elucidate various distinguishing properties. It will be seen that this sub-group provides a natural approach to the study of analytic sigmoid functions and their derivatives.

A number of integer sequences will be encountered in this note. These, along with examples of Riordan arrays, can be found in the On-Line Encyclopedia of Integer Sequences \cite{SL1, SL2}.

\section{Exponential Riordan arrays and the derivative subgroup}
Given two power series
$$g(x)=1+g_1 \frac{x}{1!} + g_2 \frac{x^2}{2!} + \cdots = \sum_{n=0}^{\infty} g_n \frac{x^n}{n!},$$
and
$$f(x)=x+f_2 \frac{x^2}{2!} + f_3 \frac{x^3}{3!}+\cdots = \sum_{n=1}^{\infty} f_n \frac{x^n}{n!}, $$
we call the matrix with general $(n,k)$-th element
$$t_{n,k}= \frac{n!}{k!} [x^n] g(x) f(x)^k$$ the exponential Riordan array associated with the pair of power series $g(x), f(x)$, and we write this matrix as $[g, f]$, or $[g(x), f(x)]$. Here, $[x^n]$ is the operator that extracts the coefficient of $x^n$ in a power series \cite{MC}. In the notation $[g(x), f(x)]$, the variable ``$x$'' is a synthetic or dummy variable, in that
$$t_{n,k}=[x^n] g(x)f(x)^k = [z^n] g(z)f(z)^k.$$

Every such matrix is lower-triangular, with $t_{n,n}=1$ for each $n \in \mathbb{N}_0$.

For the power series defined above, with their coefficients $g_n$ and $f_n$ drawn from an appropriate ring $R$, the set of matrices $[g(x), f(x)]$ forms a group, for the operation of matrix multiplication. This corresponds algebraically to the following rule.
$$ [g(x), f(x)] \cdot [u(x), v(x)] = [g(x)u(f(x)), v(f(x))].$$
The inverse element of the group of $[g(x), f(x)]$ is given by
$$[g(x), f(x)]^{-1}=\left[\frac{1}{g(\bar{f}(x))}, \bar{f}(x)\right],$$ where
$\bar{f}(x)$ is the compositional inverse of $f(x)$. This means that we have
$$f(\bar{f}(x))=\bar{f}(f(x))=x.$$
The identity element of the group is $[1, x]$.

The exponential Riordan array $[g(x), f(x)]$ has a bivariate generating function given by
$$ g(x)e^{yf(x)}=\sum_{n=0}^{\infty}\sum_{k=0}^{\infty} t_{n,k}\frac{x^n}{n!}y^k.$$
\begin{example} Regarded as an infinite lower-triangular matrix, Pascal's triangle is the element
$$[e^x, x]$$ of the group of exponential Riordan arrays (for any base ring $R$ that contains the integers $\mathbb{Z}$). This is so because we have
\begin{eqnarray*}
t_{n,k}&=& \frac{n!}{k!}[x^n] e^x x^k\\
&=& \frac{n!}{k!}[x^{n-k}] e^x\\
&=& \frac{n!}{k!}[x^{n-k}] \sum_{i=0}^{\infty} \frac{x^i}{i!}\\
&=& \frac{n!}{k!} \frac{1}{(n-k)!}\\
&=& \binom{n}{k}.\end{eqnarray*}
The inverse of this array is the array $[e^{-x}, x]$ with general element $(-1)^{n-k}\binom{n}{k}$.
\end{example}

In this note we shall be interested mainly in the \emph{derivative} sub-group of the group of exponential Riordan arrays. This subgroup consists of pairs $[f'(x), f(x)]$. That this is a sub-group follows directly from the rule for the derivative of a composed function.

\begin{example}
The element $[\cos(x), \sin(x)]$ is an element of the derivative subgroup, since $\sin'(x)=\cos(x)$. This begins
$$\left(
\begin{array}{ccccccc}
 1 & 0 & 0 & 0 & 0 & 0 & 0 \\
 0 & 1 & 0 & 0 & 0 & 0 & 0 \\
 -1 & 0 & 1 & 0 & 0 & 0 & 0 \\
 0 & -4 & 0 & 1 & 0 & 0 & 0 \\
 1 & 0 & -10 & 0 & 1 & 0 & 0 \\
 0 & 16 & 0 & -20 & 0 & 1 & 0 \\
 -1 & 0 & 91 & 0 & -35 & 0 & 1 \\
\end{array}
\right).$$
For this array, we have
\begin{eqnarray*}
t_{n,k}&=&\frac{n!}{k!}[x^n] \cos(x) \sin(x)^k\\
&=& \frac{n!}{k!}[x^n]\frac{e^{ix}+e^{-ix}}{2}\left(\frac{e^{ix}-e^{-ix}}{2i}\right)^k\\
&=& \frac{n!}{k!}[x^n]\frac{e^{ix}+e^{-ix}}{2}\frac{1}{2^k i^k}\sum_{j=0}^k\binom{k}{j} e^{ijx}(-1)^{k-j}e^{-i(k-j)x}\\
&=& \frac{n!}{k!}[x^n]\frac{e^{ix}+e^{-ix}}{2}\frac{1}{2^k i^k}\sum_{j=0}^k\binom{k}{j} (-1)^{k-j}e^{i(2j-k)x}\\
&=& \frac{n!}{k!}[x^n]\frac{1}{2^{k+1}i^k} \sum_{j=0}^k \binom{k}{j}(-1)^{k-j}\left(e^{i(2j-k+1)x}+e^{i(2j-k-1)x}\right)\\
&=& \frac{n!}{k!}\frac{1}{2^{k+1}i^k} \sum_{j=0}^k \binom{k}{j}(-1)^{k-j}\frac{i^n (2j-k+1)^n+i^n (2j-k-1)^n}{n!}\\
&=& \frac{1}{k!} \frac{1}{2^{k+1}i^k} \sum_{j=0}^k \binom{k}{j}(-1)^{k-j} \sum_{\ell=0}^n i^n \binom{n}{\ell}(2j-k)^{\ell}(1+(-1)^{n-\ell}).\end{eqnarray*}
We arrive at
$$t_{n,k}=\frac{1}{k!} \frac{1}{2^{k+1}} \sum_{j=0}^k \binom{k}{j}(-1)^{k-j} \sum_{\ell=0}^n \binom{n}{\ell}(2j-k)^{\ell}(1+(-1)^{n-\ell})(-1)^{\frac{n-k}{2}}\frac{1+(-1)^{n-k}}{2}.$$
As in this case $\cos(x)$ is an even function, and $\sin(x)$ is an odd function, this array is also an element of what is known as the ``checkerboard'' sub-group of the group of exponential Riordan arrays. The fact that $\sin(x)$ is odd is due to the fact that $\sin(-x)=-\sin(x)$, and similarly $\cos(x)$ is even since $\cos(-x)=cos(x)$, that is, $\cos(x)$ is symmetric about the line $y=0$. Such considerations will be important for some of the sigmoid functions that we will examine.

We now investigate what is known as the \emph{production matrix}\cite{ProdMat_0, ProdMat} of the exponential Riordan array $[\cos(x), \sin(x)]$. The production matrix of an exponential Riordan array $[g(x), f(x)]$ is the infinite matrix with bivariate generating function (exponential in $x$, ordinary in $y$) given by
$$ e^{xy}(Z(x)+ yA(x)),$$ \cite{Construct} where
$$A(x)=f'(\bar{f}(x)),$$
and
$$Z(x)=\frac{g'(\bar{f}(x))}{g(\bar{f}(x))}.$$
For $[g(x), f(x)]=[\cos(x), \sin(x)]$, we have
$$\bar{f}(x)=\sin^{-1}(x),$$ and we obtain
$$A(x)=\cos(\sin^{-1}(x))=\sqrt{1-x^2},\quad\text{and}\quad Z(x)=\frac{-\sin(\sin^{-1}(x))}{\cos(\sin^{-1}(x))}=\frac{-x}{\sqrt{1-x^2}}.$$
Thus the production matrix $P$ of the exponential Riordan array $[\cos(x), \sin(x)]$ has generating function
$$e^{xy}\left(\frac{-x}{\sqrt{1-x^2}}+y\sqrt{1-x^2}\right).$$ This matrix begins
$$\left(
\begin{array}{ccccccc}
 0 & 1 & 0 & 0 & 0 & 0 & 0 \\
 -1 & 0 & 1 & 0 & 0 & 0 & 0 \\
 0 & -3 & 0 & 1 & 0 & 0 & 0 \\
 -3 & 0 & -6 & 0 & 1 & 0 & 0 \\
 0 & -15 & 0 & -10 & 0 & 1 & 0 \\
 -45 & 0 & -45 & 0 & -15 & 0 & 1 \\
 0 & -315 & 0 & -105 & 0 & -21 & 0 \\
\end{array}
\right).$$
This matrix has Hessenberg form. It ``generates'' the array $[\cos(x), \sin(x)]$ in that the first row of the $n$-th power of this matrix is the $n+1$-st row of  $[\cos(x), \sin(x)]$.
The production matrix of the inverse array
$$[\cos(x), \sin(x)]^{-1}=\left[\frac{1}{\sqrt{1-x^2}}, \sin^{-1}(x)\right]$$ has generating function
$$e^{xy}(\sin(x)\sec^2(x)+y \sec(x)).$$ This matrix begins
$$\left(
\begin{array}{ccccccc}
 0 & 1 & 0 & 0 & 0 & 0 & 0 \\
 1 & 0 & 1 & 0 & 0 & 0 & 0 \\
 0 & 3 & 0 & 1 & 0 & 0 & 0 \\
 5 & 0 & 6 & 0 & 1 & 0 & 0 \\
 0 & 25 & 0 & 10 & 0 & 1 & 0 \\
 61 & 0 & 75 & 0 & 15 & 0 & 1 \\
 0 & 427 & 0 & 175 & 0 & 21 & 0 \\
\end{array}
\right).$$
\end{example}
For the Binomial matrix (Pascal's triangle) $[e^x, x]$, we can show that its production matrix has generating function $$e^{xy}(1+y).$$ This matrix begins
$$\left(
\begin{array}{ccccccc}
 1 & 1 & 0 & 0 & 0 & 0 & 0 \\
 0 & 1 & 1 & 0 & 0 & 0 & 0 \\
 0 & 0 & 1 & 1 & 0 & 0 & 0 \\
 0 & 0 & 0 & 1 & 1 & 0 & 0 \\
 0 & 0 & 0 & 0 & 1 & 1 & 0 \\
 0 & 0 & 0 & 0 & 0 & 1 & 1 \\
 0 & 0 & 0 & 0 & 0 & 0 & 1 \\
\end{array}
\right).$$
This coincides with the well-known property of the Binomial matrix
$$\binom{n+1}{k+1}=\binom{n}{k}+\binom{n}{k+1}.$$
An important property of elements of the derivative subgroup of the exponential Riordan group is the following.
\begin{proposition} The production matrix of the exponential Riordan array $[f'(x), f(x)]$ is given by
$$ U \cdot \left[\frac{1}{\bar{f}'(x)}, x\right],$$ where
$U$ is the shift matrix (with general element $\frac{n!}{(k-1)!}[n-k+1=0]$) that begins
$$\left(
\begin{array}{ccccccc}
 0 & 1 & 0 & 0 & 0 & 0 & 0 \\
 0 & 0 & 1 & 0 & 0 & 0 & 0 \\
 0 & 0 & 0 & 1 & 0 & 0 & 0 \\
 0 & 0 & 0 & 0 & 1 & 0 & 0 \\
 0 & 0 & 0 & 0 & 0 & 1 & 0 \\
 0 & 0 & 0 & 0 & 0 & 0 & 1 \\
 0 & 0 & 0 & 0 & 0 & 0 & 0 \\
\end{array}
\right).$$
\end{proposition}
The matrix $U \cdot \left[\frac{1}{\bar{f}'(x)}, x\right]$ is given by removing the first row from the Riordan array $\left[\frac{1}{\bar{f}'(x)}, x\right]$. It is thus of Hessenberg form. Note that we have used the Iverson bracket notation here: $[\mathcal{P}]$ is $1$ if the proposition $\mathcal{P}$ is true, otherwise it is $0$.
\begin{proof} The Riordan array $\left[\frac{1}{\bar{f}'(x)}, x\right]$ has bivariate generating function given by
$$ \frac{1}{\bar{f}'(x)} e^{xy}.$$
Removing the first row of the corresponding matrix corresponds to taking the derivative of this generating function with respect to $x$.
We get
$$\frac{d}{dx} \left( \frac{1}{\bar{f}'(x)} e^{xy}\right)=e^{xy}\left(\frac{d}{dx}\frac{1}{\bar{f}'(x)}+y\frac{1}{\bar{f}'(x)}\right).$$
Now $$f(\bar{f}(x))=x \Longrightarrow f'(\bar{f}(x)).\bar{f}'(x)=1$$ and hence
$$f'(\bar{f}(x))=\frac{1}{\bar{f}'(x)}.$$
Now
\begin{eqnarray*}
\frac{d}{dx} \left( \frac{1}{\bar{f}'(x)}e^{xy}\right)&=&e^{xy}\left(\frac{d}{dx}\frac{1}{\bar{f}'(x)}+y\frac{1}{\bar{f}'(x)}\right)\\
&=& e^{xy}\left(\frac{d}{dx} f'(\bar{f}(x)) + y \bar{f}'(f(x))\right)\\
&=& e^{xy} \left(f''(\bar{f}(x))\bar{f}'(x) + y \bar{f}'(f(x))\right)\\
&=& e^{xy} \left(\frac{f''(\bar{f}(x))}{f'(\bar(f)(x))}+y \bar{f}'(f(x))\right)\\
&=& e^{xy} \left(\frac{g'(\bar{f}(x))}{g(\bar{f}(x))} + y  \bar{f}'(f(x))\right),\end{eqnarray*}
where $g(x)=f'(x)$.
\end{proof}

Exponential Riordan arrays whose production matrices are of the form
$$ e^{xy} (\alpha+ \beta x + y(1+\gamma x + \delta x^2))$$ define families of orthogonal polynomials \cite{Gautschi, Chihara, Szego}.
Such production matrices are of a tri-diagonal form which begins
$$\left(
\begin{array}{ccccccc}
 \alpha  & 1 & 0 & 0 & 0 & 0 & 0 \\
 \beta  & \alpha +\gamma  & 1 & 0 & 0 & 0 & 0 \\
 0 & 2 (\beta +\delta ) & \alpha +2 \gamma  & 1 & 0 & 0 & 0 \\
 0 & 0 & 3 \beta +6 \delta  & \alpha +3 \gamma  & 1 & 0 & 0 \\
 0 & 0 & 0 & 4 \beta +12 \delta  & \alpha +4 \gamma  & 1 & 0 \\
 0 & 0 & 0 & 0 & 5 \beta +20 \delta  & \alpha +5 \gamma  & 1 \\
 0 & 0 & 0 & 0 & 0 & 6 \beta +30 \delta  & \alpha +6 \gamma  \\
\end{array}
\right).$$
This is the Jacobi matrix associated to the family of orthogonal polynomials $P_n(x)$ that satisfy the three-term recurrence relation
$$P_n(x)=(x-(\alpha+ n \gamma))P_{n-1}(x) - ((n-1)\beta+(n-1)(n-2)\delta)P_{n-2}(x),$$ with
$P_0(x)=1$, $P_1(x)=x-\alpha$.

Every family of orthogonal polynomials defines a lower-triangular matrix of coefficients $a_{n,k}$ where
$$P_n(x)=\sum_{k=0}^n a_{n,k}x^k.$$
The inverse of this coefficient array is the ``moment matrix'' of the family of orthogonal polynomials, and its first column contains the moments $m_n$ of the family. Regarded as a sequence $m_n$, we can form its Hankel transform, which is the sequence of determinants $h_n=|m_{i+j}|_{0 \le i,j \le n}$.

We finish this section by remarking that all the pairs $(f', f)$ that we shall examine are such that $f(x)$ is an odd function, and therefore $f'(x)$ is an even function. Thus the exponential Riordan array $[f'(x), f(x)]$ will be an element of the checkerboard sub-group of the Riordan group, consisting of those arrays whose first element is even and whose second element is odd. Thus each sigmoid array to follow is in the intersection of the derivative sub-group and the checkerboard sub-group.

With this review of the exponential Riordan group completed, we are now in a position to study sigmoid functions in the context of the exponential Riordan group.
\section{A Riordan analysis of sigmoid pairs}
We begin by looking at a scaled (and translated) version of the standard sigmoid function
$$ f(x)=\frac{1}{1+e^{-x}}.$$
We first translate in the $y$ direction to get
\begin{eqnarray*}
f(x)-\frac{1}{2}&=&\frac{1}{1+e^{-x}}-\frac{1}{2}\\
&=&\frac{2-1-e^{-x}}{2(1+e^{-x})}\\
&=&\frac{1}{2} \frac{1-e^{-x}}{1+e^{-x}}\\
&=& \frac{1}{2} \frac{e^{\frac{x}{2}}-e^{\frac{-x}{2}}}{e^{\frac{x}{2}}+e^{\frac{-x}{2}}}\\
&=& \frac{1}{2} \tanh\left(\frac{x}{2}\right).\end{eqnarray*}
Now $\frac{1}{2} \tanh\left(\frac{x}{2}\right)=\sum_{n=0}^{\infty}a_n \frac{x^n}{n!}$, where the sequence $a_n$ begins
$$0, \frac{1}{4}, 0, - \frac{1}{8}, 0, \frac{1}{8}, 0, - \frac{17}{16}, 0, \frac{31}{4}, 0,\ldots.$$
We prefer to work with integer valued coefficients, for simplicity, so we rescale again to study the function
$$\tilde{f}(x)=\frac{2}{1+e^{-2x}}.$$
We then have
$$\tilde{f}(x)-1=\frac{2}{1+e^{-2x}}-1=\tanh(x).$$

\begin{center}
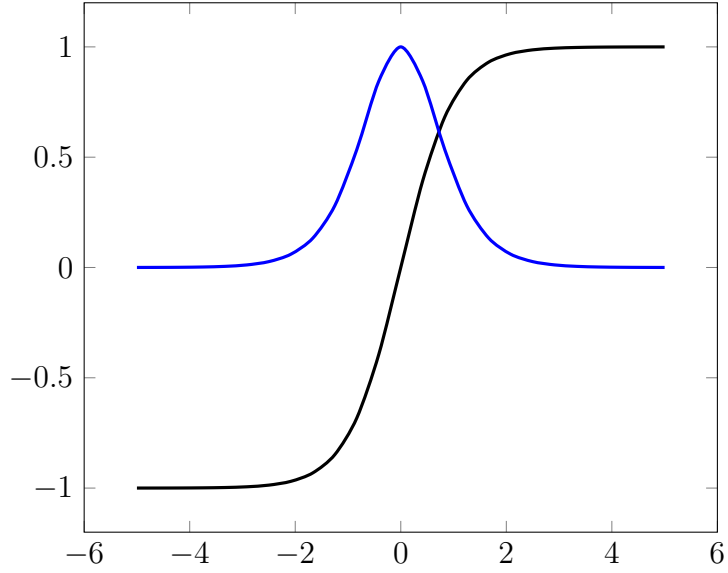
\begin{figure}
\begin{center}
\begin{tikzpicture}
\begin{axis} [width=10cm]
\addplot[smooth, very thick, color=black] {(exp(x)-exp(-x))/(exp(x)+exp(-x))};
\addplot[smooth, very thick, color=blue] {4*exp(2*x)/(1+exp(2*x))^2};
\end{axis}
\end{tikzpicture}
\end{center}
\caption{The sigmoid pair $\sech^2(x), \tanh(x)$}\label{sigmoid_2}
\end{figure}
\end{center}
Thus we commence our study of sigmoid functions with the pair $(\tanh'(x), \tanh(x))=(\sech^2(x), \tanh(x))$. See Figure \ref{sigmoid_2}.

The associated exponential Riordan array, which we shall call a \emph{sigmoid array}, is
$$[\sech^2(x), \tanh(x)] = \left[\frac{1}{1-x^2}, \tanh^{-1}(x)\right]^{-1}.$$
The production array of $[\sech^2(x), \tanh(x)]$ has generating function
$$e^{xy}\left(-2x+y(1-x^2)\right).$$
This implies that the inverse of sigmoid array $[\sech^2(x), \tanh(x)]$ is the coefficient array of the family of formally orthogonal polynomials that satisfy the three-term recurrence
$$P_n(x)=xP_{n-1}(x)-n(1-n)P_{n-2}(x).$$
These polynomials begin
$$1, x, x^2 + 2, x(x^2 + 8), x^4 + 20x^2 + 24, x(x^4 + 40x^2 + 184), x^6 + 70x^4 + 784x^2 + 720,\ldots.$$
The production matrix of the sigmoid array $[\sech^2(x), \tanh(x)]$  is given by
$$ P=U.[1-x^2, x].$$
We deduce that the Hankel transform \cite{Hankel_Riordan, Kratt, KrattC, Layman} of the expansion of $f'(x)=\sech^2(x)$ is given by
$$h_n=\prod_{k=0}^n ((k+2)(1-(k+2)))^{n-k},$$ which begins
$$1, -2, -24, 3456, 9953280, -859963392000, -3120635156889600000,\ldots.$$
We also have that the Hankel transform of the expansion of sigmoid function $\tanh(x)$ is given by
$$\prod_{k=0}^n k!^2 (-1)^{\frac{n+1}{2}}\frac{1-(-1)^n}{2}.$$ This sequence begins
$$0, -1, 0, 144, 0, -1194393600, 0, 15728001190723584000000,\ldots$$

\begin{center}
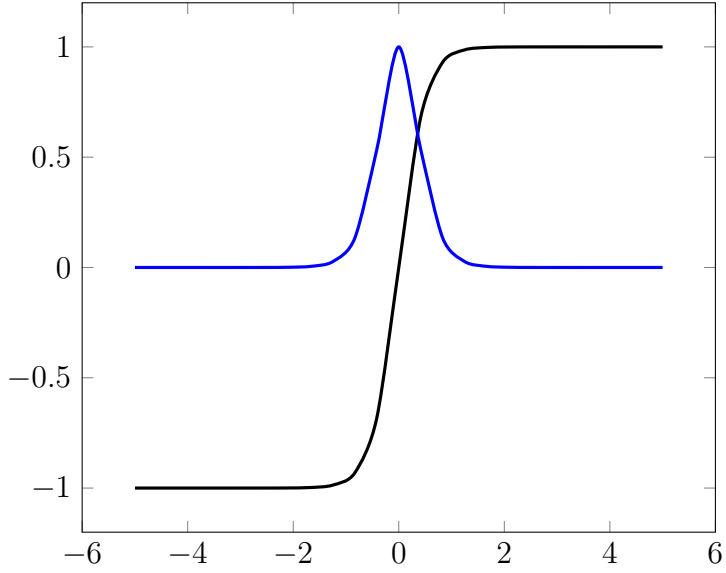
\begin{figure}
\begin{center}
\begin{tikzpicture}
\begin{axis} [width=10cm]
\addplot[smooth, very thick, color=black] {(1/2)(exp(2*x)-exp(-2*x))/(exp(2*x)+exp(-2*x))};
\addplot[smooth, very thick, color=blue] {4*exp(4*x)/(1+exp(4*x))^2};
\end{axis}
\end{tikzpicture}
\end{center}
\caption{The sigmoid pair $\sech^2(2x),\frac{1}{2} \tanh(2x)$}\label{sech2}
\end{figure}
\end{center}
A related sigmoid array is the array
$$[\sech^2(2x), \frac{1}{2} \tanh(2x)]=\left[\frac{1}{1-4x^2}, \frac{1}{2} \tanh^{-1}(2x)\right]^{-1},$$ based on the sigmoid function $\frac{1}{2} \tanh(2x)$. See Figure \ref{sech2}.

This sigmoid array has its production matrix generated by
$$e^{xy}(-8x+y(1-4x^2)),$$ and hence its inverse array is the coefficient array of a family of formally orthogonal polynomials $Q_n(x)$ that satisfy the three-term recurrence
$$Q_n(x)=xQ_{n-1}(x)-4n(1-n)Q_{n-2}(x).$$
The production matrix of this sigmoid array is given by
$$P=U\cdot [1-4x^2, x].$$

The next sigmoid array that we shall study is given by
$$\left[\frac{1}{1+x^2}, \tan^{-1}(x)\right]=[1+\tan^2(x), \tan(x)]^{-1}=[\sec^2(x), \tan(x)]^{-1},$$ based on the sigmoid function $\tan^{-1}(x)$. See Figure \ref{atan}.
\begin{center}
\begin{figure}
\begin{center}
\begin{tikzpicture}
\begin{axis} [width=10cm]
\addplot[smooth, very thick, color=black] {rad(atan(x))};
\addplot[smooth, very thick, color=blue] {1/(1+x^2)};
\end{axis}
\end{tikzpicture}
\end{center}
\caption{The sigmoid pair $\frac{1}{1+x^2},\tan^{-1}(x)$}\label{atan}
\end{figure}
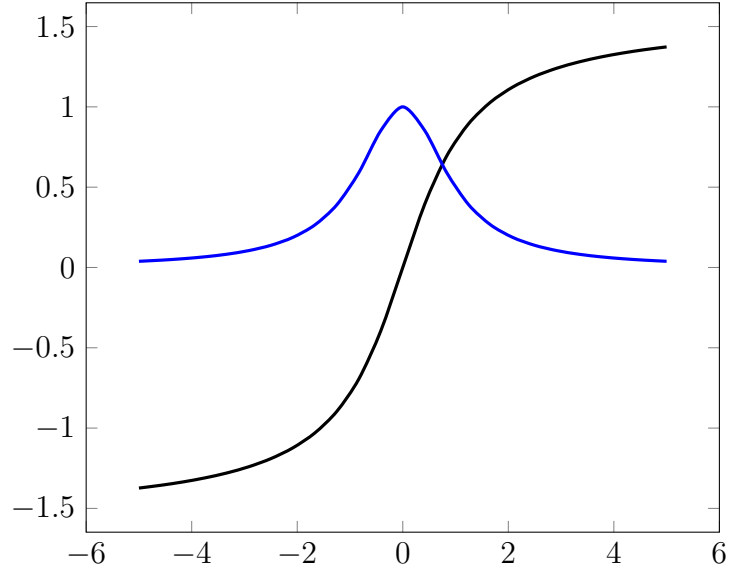
\end{center}
We find that the generating function of the inverse array $[\sec^2(x), \tan(x)]$ is given by
$$e^{xy}\left(2x+y(1+x^2)\right),$$  and thus the sigmoid array $\left[\frac{1}{1+x^2}, \tan^{-1}(x)\right]$ is itself the coefficient array of the family of orthogonal polynomials $P_n(x)$ that satisfy the three-term recurrence
$$P_n(x)=xP_{n-1}(x)-n(n-1)P_{n-2}(x),$$ with
$P_0(x)=1$, $P_1(x)=x$. These polynomials begin
$$1, x, x^2 - 2, x(x^2 - 8), x^4 - 20x^2 + 24, x(x^4 - 40x^2 + 184), x^6 - 70x^4 + 784x^2 - 720,\ldots.$$
The moments of this family of polynomials, which begin
$$1, 0, 2, 0, 16, 0, 272, 0, 7936, 0, 353792, 0, 22368256,0,\ldots$$ therefore have their generating function given by $\sec^2(x)$, with
Hankel transform given by $$\prod_{k=0}^n ((k+1)(k+2))^{n-k}.$$
The production matrix of this sigmoid array is given by
$$U\cdot [\cos^2(x), x].$$

We next consider the sigmoid array
$$\left[\frac{1}{(1+x^2)^{\frac{3}{2}}}, \frac{x}{\sqrt{1+x^2}}\right]=\left[\frac{1}{(1-x^2)^{\frac{3}{2}}}, \frac{x}{\sqrt{1-x^2}}\right]^{-1},$$ based on the sigmoid function $\frac{x}{\sqrt{1+x^2}}$. See Figure \ref{power}.

\begin{center}
\begin{figure}
\begin{center}
\begin{tikzpicture}
\begin{axis} [width=10cm]
\addplot[smooth, very thick, color=black] {1/(1+x^2)^(3/2)};
\addplot[smooth, very thick, color=blue] {x/(1+x^2)^(1/2)};
\end{axis}
\end{tikzpicture}
\end{center}
\caption{The sigmoid pair $\frac{1}{(1+x^2)^{3/2}}, \frac{x}{\sqrt{1+x^2}}$}\label{power}
\end{figure}
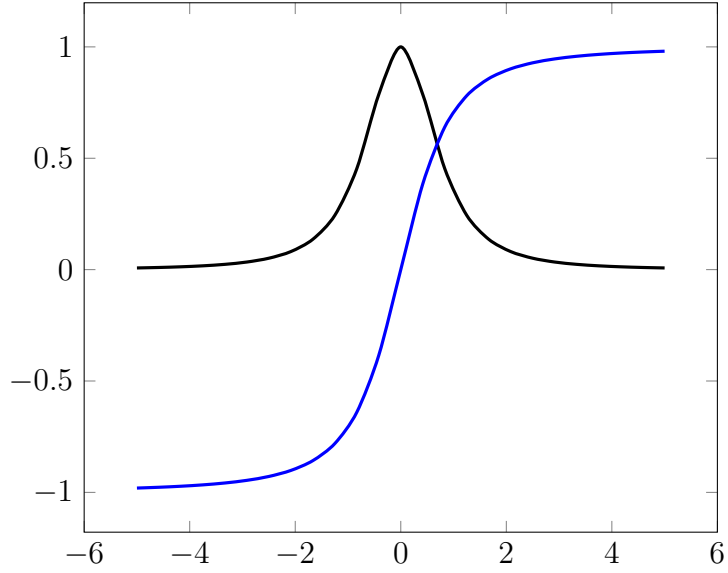
\end{center}
The generating function of the production matrix of the sigmoid array is
$$e^{xy}(-3x \sqrt{1-x^2}+y(1-x^2)^{\frac{3}{2}}).$$
Similarly, the generating function of the inverse array is given by
$$e^{xy}(3x \sqrt{1+x^2}+y(1+x^2)^{\frac{3}{2}}).$$
Clearly, neither array is the coefficient array of a family of orthogonal polynomials. Thus this sigmoid pair is qualitatively different from the previous examples.
The production matrix of this sigmoid array is given by
$$U \cdot \left[ (1-x^2)^{\frac{3}{2}}, x\right].$$

A related sigmoid array (one of an obvious family) is given by
$$\left[\frac{1}{(1+x^4)^{\frac{5}{4}}}, \frac{x}{(1+x^4)^{\frac{1}{4}}}\right]=\left[\frac{1}{(1-x^4)^{\frac{5}{4}}}, \frac{x}{(1-x^4)^{\frac{1}{4}}}\right]^{-1}.$$
We find that the production matrix of the sigmoid array has generating function
$$e^{xy}(-5x^3(1-x^4)^{\frac{1}{4}}+y(1-x^4)^{\frac{5}{4}}),$$ with an analogous result for the inverse array.

The production matrix of this sigmoid array is given by
$$U \cdot \left[(1-x^4)^{\frac{5}{4}},x\right].$$

We next consider the sigmoid array given by
$$[\sech(x), \gd(x)]$$ where the Gudermanian function $\gd(x)$ is defined by
$$\gd(x)=\tan^{-1}(\sinh(x))=\sin^{-1}(\tanh(x))=\cos^{-1}(\sech(x)).$$ See Figure \ref{gud}.
\begin{center}
\begin{figure}
\begin{center}
\begin{tikzpicture}
\begin{axis} [width=10cm]
\addplot[smooth, very thick, color=black] {rad(asin(tanh(x))};
\addplot[smooth, very thick, color=blue] {2*exp(x)/(exp(2*x)+1)};
 \end{axis}
\end{tikzpicture}
\end{center}
\caption{The sigmoid pair $\sec(x), \gd(x)$}\label{gud}
\end{figure}
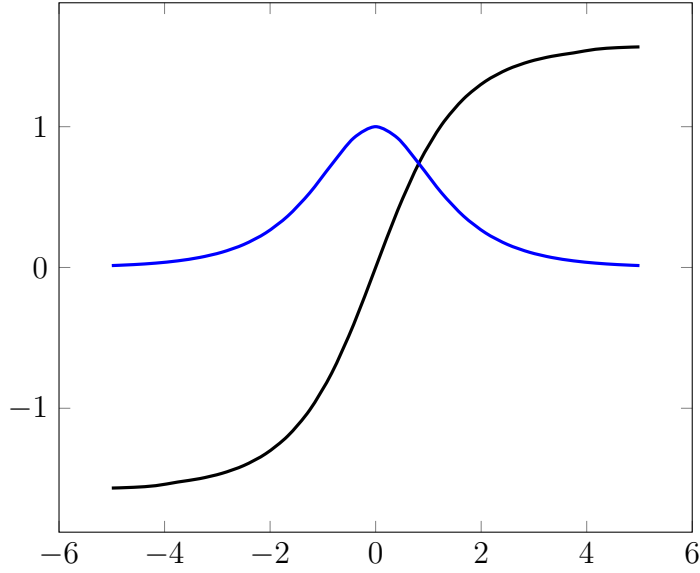
\end{center}
The inverse of this sigmoid array may be expressed in many ways. Some are
$$[\sech(x), \gd(x)]^{-1}=[\sec(x), \sinh^{-1}(\tan(x))]=[\sec(x), \log(\sec(x)+\tan(x))],$$
and
$$[\sech(x), \gd(x)]^{-1}=[\sec(x), \int_0^x \sec(t)\,dt].$$
Here, $\gd'(x)=\sech(x)$ and $\gd^{-1}(x)=\sech^{-1}(\cos(x))$. Using these expressions we can show that the generating function of the production matrix of the sigmoid array $[\sech(x), \gd(x)]$ is given by
$$e^{xy}(-\sin(x)+y \cos(x)).$$ In a similar manner, we can show that the generating function of the inverse sigmoid array is given by
$$e^{xy} (\sinh(x)+y \cosh(x)).$$
Clearly, neither of these represent the production array of a coefficient array of orthogonal polynomials.
However, we note that the exponential Riordan array $[\sech(x), \tanh(x)]$ is the moment array (inverse of the coefficient array) of the formally orthogonal polynomial family given by
$$P_n(x)=xP_{n-1}(x)+(n-1)^2P_{n-2}(x),$$ with
$P_0(x)=1$, $P_1(x)=x$. We also have
\begin{eqnarray*}
[\sech(x), \gd(x)]^{-1}\cdot [\sech(x),\tan(x)]&=&\left[\frac{1}{\cos(x)}, \sinh^{-1}(\tan(x))\right]\cdot[\sech(x), \tanh(x)]\\
&=&\left[\frac{1}{\cos(x)} \sech(\sinh^{-1}(\tan(x))), \tanh(\sinh^{-1}(\tan(x)))\right]\\
&=&[1, \sin(x)].\end{eqnarray*}
We therefore obtain the relation
$$[\sech(x), \gd(x)]=[\sech(x), \tanh(x)] \cdot [1, \sin^{-1}(x)],$$ relating the sigmoid matrix
$[\sech(x), \gd(x)]$ to the moment array of the above family of orthogonal polynomials.
Equivalently, we have
$$[\sech(x), \gd(x)]^{-1}=[1, \sin(x)]\cdot [\sech(x), \tanh(x)]^{-1}$$ which shows that the coefficient array
$[\sech(x), \tanh(x)]^{-1}$ of the above family of orthogonal polynomials is related to the inverse of the sigmoid array by the invertible transformation $[1, \sin(x)]$.

The production matrix of the sigmoid array $[\sech(x), \gd(x)]$ is given by
$$U \cdot [\cos(x), x].$$

We next consider the sigmoid array
$$\left[e^{-x^2}, \frac{\sqrt{\pi}}{2} \erf(x)\right],$$ based on the sigmoid function $\frac{\sqrt{\pi}}{2} \erf(x)$. See Figure \ref{erf}.

\makeatletter
\pgfmathdeclarefunction{erf}{1}{%
  \begingroup
    \pgfmathparse{#1 > 0 ? 1 : -1}%
    \edef\sign{\pgfmathresult}%
    \pgfmathparse{abs(#1)}%
    \edef\x{\pgfmathresult}%
    \pgfmathparse{1/(1+0.3275911*\x)}%
    \edef\t{\pgfmathresult}%
    \pgfmathparse{%
      1 - (((((1.061405429*\t -1.453152027)*\t) + 1.421413741)*\t
      -0.284496736)*\t + 0.254829592)*\t*exp(-(\x*\x))}%
    \edef\y{\pgfmathresult}%
    \pgfmathparse{(\sign)*\y}%
    \pgfmath@smuggleone\pgfmathresult%
  \endgroup
}
\makeatother
\begin{center}
\begin{figure}
\begin{center}
\begin{tikzpicture}[yscale = 3]
  \draw[thick, ->] (-5,0) -- node[at end,below] {$x$}(5,0);
  \draw[thick, ->] (0,-1) -- node[below left] {$0$} node[at end,
  left] {$erf(x)$} (0,1);
  \draw[black, very thick] plot[domain=-5:5,samples=200] (\x, {(0.8862269254)*erf(\x)});
  \draw[blue, very thick] plot[domain=-5:5,samples=200] (\x,{exp(-(\x)^2)});
\end{tikzpicture}
\caption{The sigmoid pair $e^{-x^2}, \frac{\sqrt{\pi}}{2} \erf(x)$}\label{erf}
\end{center}
\end{figure}
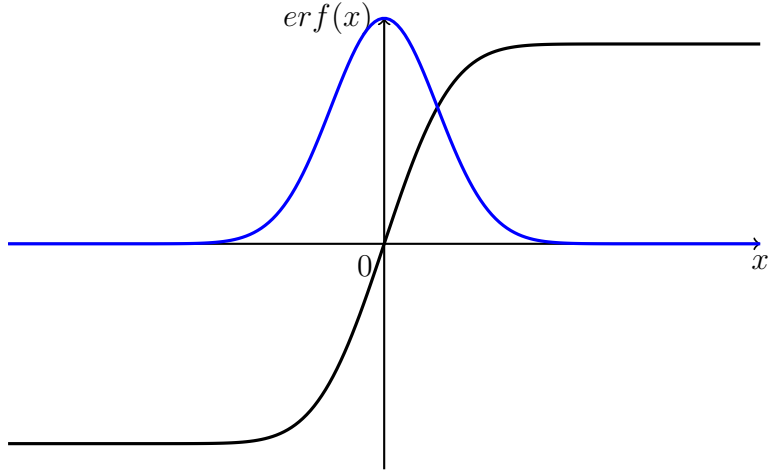
\end{center}

This array begins
$$\left(
\begin{array}{ccccccc}
 1 & 0 & 0 & 0 & 0 & 0 & 0 \\
 0 & 1 & 0 & 0 & 0 & 0 & 0 \\
 -2 & 0 & 1 & 0 & 0 & 0 & 0 \\
 0 & -8 & 0 & 1 & 0 & 0 & 0 \\
 12 & 0 & -20 & 0 & 1 & 0 & 0 \\
 0 & 112 & 0 & -40 & 0 & 1 & 0 \\
 -120 & 0 & 532 & 0 & -70 & 0 & 1 \\
\end{array}
\right),$$ and it has a production matrix that begins
$$\left(
\begin{array}{ccccccc}
 0 & 1 & 0 & 0 & 0 & 0 & 0 \\
 -2 & 0 & 1 & 0 & 0 & 0 & 0 \\
 0 & -6 & 0 & 1 & 0 & 0 & 0 \\
 -4 & 0 & -12 & 0 & 1 & 0 & 0 \\
 0 & -20 & 0 & -20 & 0 & 1 & 0 \\
 -56 & 0 & -60 & 0 & -30 & 0 & 1 \\
 0 & -392 & 0 & -140 & 0 & -42 & 0 \\
\end{array}
\right).$$
Clearly, this is not the production matrix of a moment matrix for a family of orthogonal polynomials (it is not tri-diagonal). However, we can again apply an invertible transform to this matrix to obtain the moment matrix of a family of orthogonal polynomials.

Thus we consider the exponential Riordan array $\left[e^{-x^2}, x\right]$. This matrix has a production matrix with generating function
$$e^{xy}(-2x+y),$$ and hence it is the moment matrix for the family of orthogonal polynomials given by
$$P_n(x)=xP_{n-1}(x)+2(n-1)P_{n-2}(x).$$
These polynomials begin
$$1, x, x^2 + 2, x^3 + 6x, x^4 + 12x^2 + 12, x^5 + 20x^3 + 60x, x^6 + 30x^4 + 180x^2 + 120,\ldots.$$
We have
$$\left[e^{-x^2}, \frac{\sqrt{\pi}}{2} \erf(x)\right]=\left[e^{-x^2},x\right]\cdot \left[1, \frac{\sqrt{\pi}}{2} \erf(x)\right].$$
In this case, the lack of a closed form for the inverse of $\erf(x)$ limits the amount of analysis that can be carried out.

Our final example differs from the previous one in that the sigmoid function we consider is not an odd function - we do not have $g(-x)=-g(x)$.  It finds application in many areas of applied mathematics including actuarial science and population dynamics. This is the Gompertz function \cite{Rad}
$$ g(x)=e^{1-e^{-x}}-1$$ which along with its derivative
$$ g'(x)= e^{1-x-e^{-x}}$$ defines the sigmoid array (which we shall call the \emph{Gompertz sigmoid array})
$$\left[e^{1-x-e^{-x}}, e^{1-e^{-x}}-1\right]=\left[\frac{1}{(1+x)(1-\log(1+x))}, \log\left(\frac{1}{1-\log(1+x)}\right)\right]^{-1}$$ which begins
$$\left(
\begin{array}{ccccccc}
 1 & 0 & 0 & 0 & 0 & 0 & 0 \\
 0 & -1 & 0 & 0 & 0 & 0 & 0 \\
 -1 & 0 & 1 & 0 & 0 & 0 & 0 \\
 1 & 4 & 0 & -1 & 0 & 0 & 0 \\
 2 & -3 & -10 & 0 & 1 & 0 & 0 \\
 -9 & -22 & 5 & 20 & 0 & -1 & 0 \\
 9 & 50 & 112 & -5 & -35 & 0 & 1 \\
\end{array}
\right).$$ 
\begin{center}
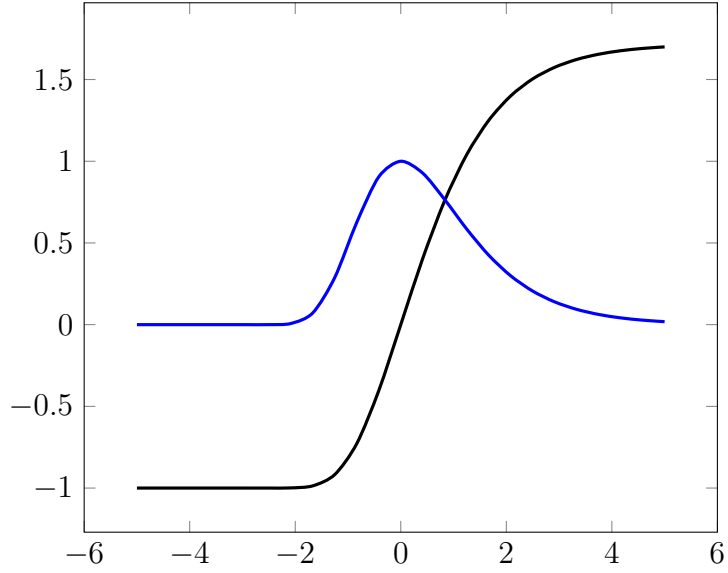
\begin{figure}
\begin{center}
\begin{tikzpicture}
\begin{axis} [width=10cm]
\addplot[smooth, very thick, color=black] {exp(1-exp(-x))-1};
\addplot[smooth, very thick, color=blue] {exp(1-x-exp(-x))};
 \end{axis}
\end{tikzpicture}
\end{center}
\caption{The Gompertz sigmoid pair $e^{1-x-e^{-x}}, e^{1-e^{-x}}-1$}\label{Gompertz}
\end{figure}
\end{center}

See Figure \ref{Gompertz}.

The production matrix of this array is given by
$$ U\cdot [(1+x)(1-\log(1+x)), x].$$
The elements of the first column of the Gompertz sigmoid array, which have exponential generating function $g'(x)=e^{1-x-e^{-x}}$ are the moments
for the family of formally orthogonal polynomials given by
$$P_n(x)=(x+(n-1))P_{n-1}(x)+(n-1)P_{n-2}(x),$$ with
$P_0(x)=1$, $P_1(x)=1$.
This means that the moment sequence
$$1,0,-1,12,-9,9,\ldots$$ has an ordinary generating function given by the continued fraction \cite{Wall}
$$\cfrac{1}{1+\cfrac{x^2}{1+x+\cfrac{2x^2}{1+2x+\cfrac{3x^2}{1+3x+\cfrac{4x^2}{1+\cdots}}}}}.$$
The moment matrix of these polynomials is given by the exponential Riordan array
$$\left[e^{1-x-e^{-x}}, 1-e^{-x}\right],$$  and we have the relation
$$\left[e^{1-x-e^{-x}}, e^{1-e^{-x}}-1\right]=\left[e^{1-x-e^{-x}}, 1-e^{-x}\right]\cdot \left[1, e^x-1\right].$$
The right-most exponential array is the array of Stirling numbers $S2(n,k)$ of the second kind, which begins
$$\left(
\begin{array}{ccccccc}
 1 & 0 & 0 & 0 & 0 & 0 & 0 \\
 0 & 1 & 0 & 0 & 0 & 0 & 0 \\
 0 & 1 & 1 & 0 & 0 & 0 & 0 \\
 0 & 1 & 3 & 1 & 0 & 0 & 0 \\
 0 & 1 & 7 & 6 & 1 & 0 & 0 \\
 0 & 1 & 15 & 25 & 10 & 1 & 0 \\
 0 & 1 & 31 & 90 & 65 & 15 & 1 \\
\end{array}
\right).$$
Thus we have
$$g_{n,k}=\sum_{j=0}^n m_{n,j}S2_{j,k}$$ where $m_{n,k}$ are the moment matrix elements for the orthogonal polynomials above, and where $g_{n,k}$ is the general element of the Gombertz sigmoid array.

In fact, it is possible to define $g_{n,k}$ solely in terms of the Stirling numbers of the second kind. Thus we have
$$g_{n,k}=\sum_{j=0}^n S2(n+1,j+1)(-1)^{n-j}S2(j+1,k+1),$$ which follows from the equality of exponential Riordan arrays
$$ \left[e^{1-x-e^{-x}}, e^{1-e^{-x}}-1\right]=\left[e^{-x}, 1-e^{-x}\right]\cdot \left[ e^x, e^x-1\right].$$
In particular, the elements in the expansion of the derivative of the Gompertz sigmoid function are given by
$$g_{n,0}= \sum_{j=0}^n S2(n+1,j+1)(-1)^{n-j}S2(j+1,1)=\sum_{j=0}^n S2(n+1,j+1)(-1)^{n-j}.$$

\section{Parametric visualizations}

We can visualize the exponential Riordan arrays $[f'(x), f(x)]$ that we have studied by providing parametric plots of the pair $(f'(t), f(t))$. For the example of $[\cos(x), \sin(x)]$, we get the unit circle. We provide visualizations of the sigmoid arrays that we have considered.

\begin{figure}
\begin{center}
\begin{tikzpicture}
    \begin{axis}[
            xmin=-2,xmax=2,
            ymin=-2,ymax=2,
            grid=both,
            ]
            \addplot [very thick, domain=0:360,samples=100]({cos(x)},{sin(x)});
    \end{axis}
\end{tikzpicture}
\caption{Parametric plot of $t \mapsto (\cos(t), \sin(t))$}
\end{center}
\end{figure}
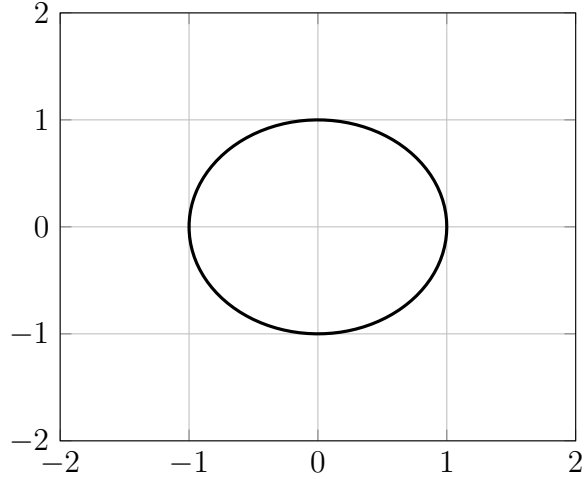

\begin{figure}
\begin{center}
\begin{tikzpicture}
    \begin{axis}[
            xmin=-2,xmax=2,
            ymin=-2,ymax=2,
            grid=both,
            ]
            \addplot [very thick, domain=-4:4,samples=100]({4*exp(2*x)/(1+exp(2*x))^2},{(1-exp(-2*x))/(1+exp(-2*x))});
    \end{axis}
\end{tikzpicture}
\caption{Parametric plot of $t \mapsto (\sech^2(t), \tanh(t))$}
\end{center}
\end{figure}
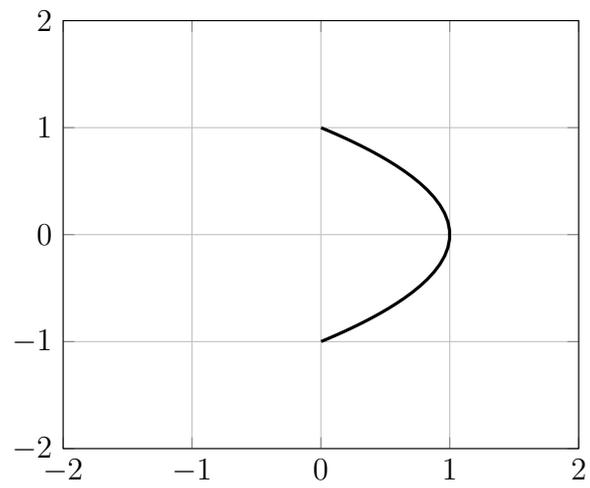

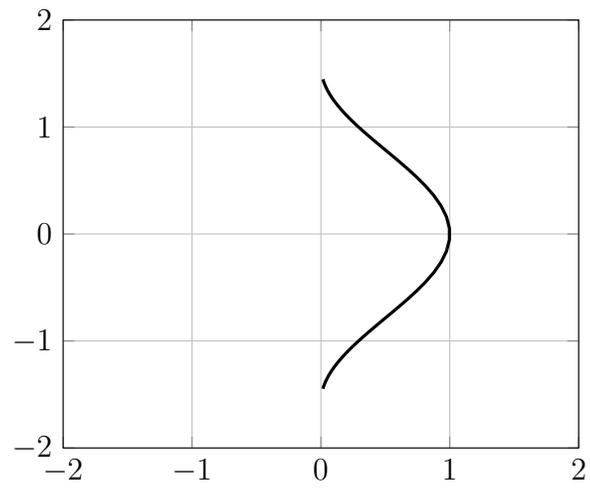
\begin{figure}
\begin{center}
\begin{tikzpicture}
    \begin{axis}[
            xmin=-2,xmax=2,
            ymin=-2,ymax=2,
            grid=both,
            ]
            \addplot [very thick, domain=-8:8,samples=150]({1/(1+x^2)},{rad(atan(x))});
    \end{axis}
\end{tikzpicture}
\caption{Parametric plot of $t \mapsto \left(\frac{1}{1+t^2}, \tan^{-1}(t)\right)$}
\end{center}
\end{figure}

\begin{figure}
\begin{center}
\begin{tikzpicture}
    \begin{axis}[
            xmin=-2,xmax=2,
            ymin=-2,ymax=2,
            grid=both,
            ]
            \addplot [very thick, domain=-4:4,samples=100]({1/(1+x^2)^(3/2)},{x/(1+x^2)^(1/2)});
    \end{axis}
\end{tikzpicture}
\caption{Parametric plot of $t \mapsto \left(\frac{1}{(1+t^2)^{3/2}}, \frac{t}{(1+t^2)^{1/2}}\right)$}
\end{center}
\end{figure}
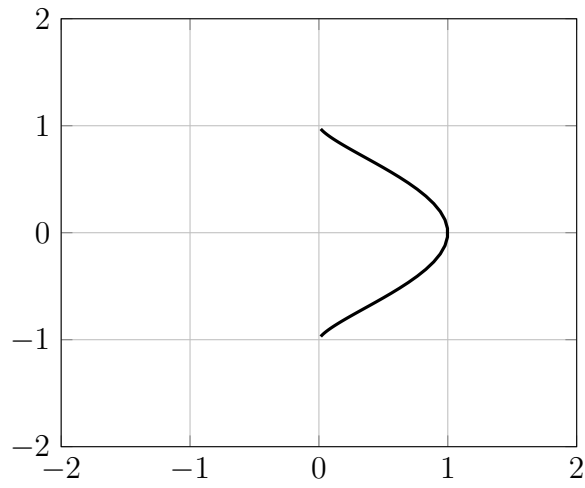

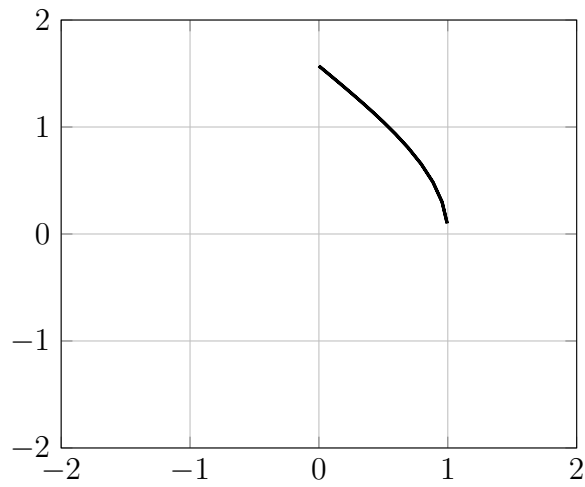
\begin{figure}
\begin{center}
\begin{tikzpicture}
    \begin{axis}[
            xmin=-2,xmax=2,
            ymin=-2,ymax=2,
            grid=both,
            ]
            \addplot [very thick, domain=-10:10,samples=100]({2*exp(x)/(1+exp(2*x))},{rad(acos(2*exp(x)/(1+exp(2*x))))});
    \end{axis}
\end{tikzpicture}
\caption{Parametric plot of $t \mapsto (\sech(t), \gd(t))$}
\end{center}
\end{figure}

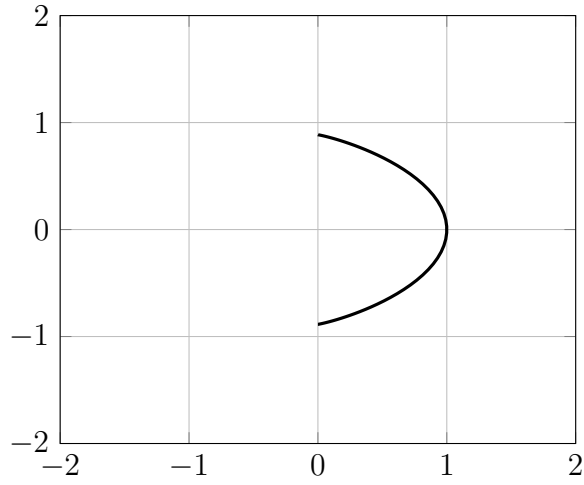
\begin{figure}
\begin{center}
\begin{tikzpicture}
    \begin{axis}[
            xmin=-2,xmax=2,
            ymin=-2,ymax=2,
            grid=both,
            ]
            \addplot [very thick, domain=-3:3,samples=100]({exp(-x^2)},{(0.8862269254)*erf(x)});
    \end{axis}
\end{tikzpicture}
\caption{Parametric plot of $t \mapsto \left(e^{-t^2}, \frac{\sqrt{\pi}}{2} \erf(t)\right)$}
\end{center}
\end{figure}

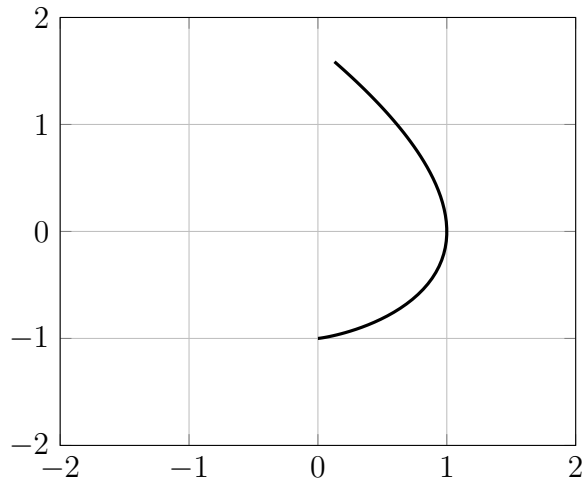
\begin{figure}
\begin{center}
\begin{tikzpicture}
    \begin{axis}[
            xmin=-2,xmax=2,
            ymin=-2,ymax=2,
            grid=both,
            ]
            \addplot [very thick, domain=-3:3,samples=100]({exp(1-x-exp(-x))},{exp(1-exp(-x))-1});
    \end{axis}
\end{tikzpicture}
\caption{Parametric plot of $t \mapsto \left(e^{1-t-e^{-t}}, e^{1-e^{-t}}-1\right)$}
\end{center}
\end{figure}

\newpage
\section{Polynomial families}
\begin{figure}
\begin{center}
\begin{tikzpicture}
\begin{axis} [width=10cm]
\addplot[smooth, very thick, color=black] {1};
\addplot[smooth, very thick, color=blue] {x};
\addplot[smooth, very thick, color=red] {x^2-3};
\addplot[smooth, very thick, color=black] {x^3-12*x};
\addplot[smooth, very thick, color=blue] {x^4-30*x^2+45};
\end{axis}
\end{tikzpicture}
\end{center}
\caption{First five polynomials defined by the sigmoid array $\left[\frac{1}{(1+x^2)^{\frac{3}{2}}}, \frac{x}{(1+x^2)^{\frac{1}{2}}}\right]$}\label{polys}
\end{figure}
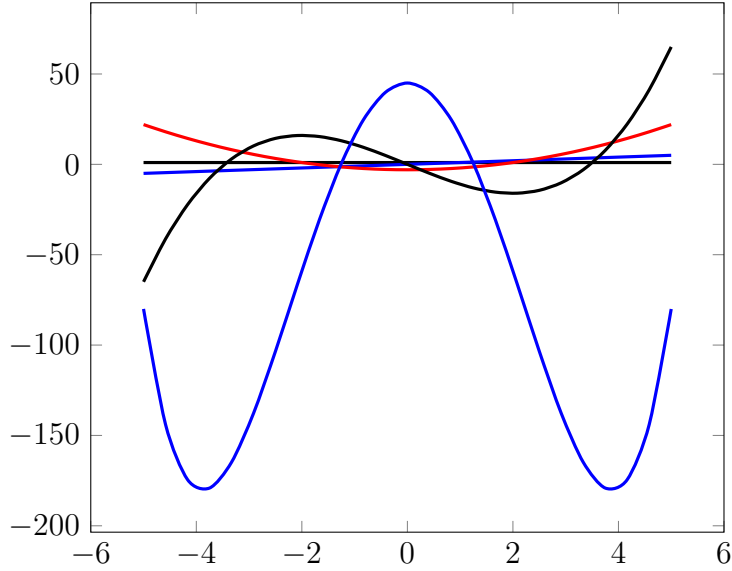

We have seen that for some sigmoid arrays, either they or their inverses define a family of orthogonal polynomials. Due to the lower-triangular nature of Riordan arrays, a family of polynomials is always produced, when the Riordan array multiplies the vector $(1,x, x^2, \ldots)^T$.
For instance, the sigmoid array $\left(\frac{1}{(1+x^2)^{\frac{3}{2}}}, \frac{x}{\sqrt{1+x^2}}\right)$, which begins
$$\left(
\begin{array}{ccccccc}
 1 & 0 & 0 & 0 & 0 & 0 & 0 \\
 0 & 1 & 0 & 0 & 0 & 0 & 0 \\
 -3 & 0 & 1 & 0 & 0 & 0 & 0 \\
 0 & -12 & 0 & 1 & 0 & 0 & 0 \\
 45 & 0 & -30 & 0 & 1 & 0 & 0 \\
 0 & 360 & 0 & -60 & 0 & 1 & 0 \\
 -1575 & 0 & 1575 & 0 & -105 & 0 & 1 \\
\end{array}
\right),$$
produces the family of polynomials that begins
$$1, x, x^2 - 3, x^3 - 12x, x^4 - 30x^2 + 45, x^5 - 60x^3 + 360x, x^6 - 105x^4 + 1575x^2 - 1575,\ldots.$$
See Figure \ref{polys}. Further studies may encompass a detailed examination of the properties of these polynomials.

\section{Conclusions} We have studied a range of sigmoid functions within the novel setting of exponential Riordan arrays. This has proven fruitful, bringing into play links to families of orthogonal polynomials in some cases, and suggesting that the new sigmoid arrays may define interesting families of polynomials in their own right. Both similarities and differences have been found in the different functions studied. The links between the Gompertz function and the Stirling numbers of the second kind, studied by Rzadkowski et  al \cite{Rad}, for instance, are made clear within the context exponential Riordan arrays. It is hoped that further investigation using these methods will lead to other significant results.

\bigskip
\hrule
\bigskip
\noindent 2010 {\it Mathematics Subject Classification}: Primary
33B10; Secondary 33B20, 42C05, 11C20, 15B36, 05A15
\noindent \emph{Keywords:} sigmoid function, special functions, exponential Riordan array, orthogonal polynomials, Hankel transform, generating functions.

\end{document}